\RequirePackage[leqno]{amsmath}
\documentclass[leqno]{amsart}
\usepackage{pifont}
\usepackage{amsmath}
\usepackage{color}
\usepackage{amsfonts}
\usepackage{enumitem}
\setlist[itemize]{topsep=0pt,after=\vspace{1.5\baselineskip}}
 \usepackage[colorlinks=true]{hyperref}
\usepackage{empheq}
\usepackage[T1]{fontenc}
\newtheorem*{state}{Assumptions}
\def\R{\mathbb R}  
\newtheorem{theorem}{Theorem}[section]

\newtheorem{lemma}[theorem]{Lemma}

\newtheorem{remark}{Remark}
\def\nx{{\bf x}}

\def\io{\int_\Omega} \def\iob{\int_{\partial\Omega}}

\title[Properties of solutions to porous medium problems] 
      {Properties of solutions to porous medium problems with different sources and boundary conditions}
\author[T. Li, N. Pintus and G. Viglialoro]{}
\keywords{Nonlinear parabolic problems, porous medium equations, global existence, blow-up, gradient nonlinearities.\\
\textit{$^\dagger$Corresponding author}: giuseppe.viglialoro@unica.it}

\begin{document}
\maketitle

\centerline{\scshape Tongxing Li$^{1,2}$, Nicola Pintus$^{3}$ \and Giuseppe Viglialoro$^{4,{\dagger}}$}
\medskip
{
 \footnotesize
  \centerline{$^1$LinDa Institute of Shandong Provincial }
  \centerline{Key Laboratory of Network Based Intelligent Computing}
  \centerline{Linyi University }
  \centerline{Linyi, Shandong, 276005 (P. R. China)}
 
 \centerline{$^2$School of Control Science and Engineering}
 \centerline{Shandong University}
 \centerline{Jinan, Shandong, 250061 (P. R. China)}
 \medskip   
 \centerline{$^3$Osservatorio Astronomico di Cagliari}
 \centerline{Cagliari, Via della Scienza 5, 09047 (Italy)}
 \medskip
 \centerline{$^4$Dipartimento di Matematica e Informatica}
  \centerline{Universit\`{a} di Cagliari}
  \centerline{Cagliari, V. le Merello 92, 09123 (Italy)}
}

\bigskip
\begin{abstract}
In this paper we study nonnegative and classical solutions $u=u(\nx,t)$ to porous medium problems of the type 
\begin{equation}\label{ProblemAbstract}
\tag{$\Diamond$}
\begin{cases}
u_t=\Delta u^m + g(u,|\nabla u|) & {\bf x} \in \Omega, t\in I,\\
u({\bf x},0)=u_0({\bf x})&{\bf x} \in  \Omega,\\
\end{cases}
\end{equation}
where $\Omega$ is a bounded and smooth domain of $\R^N$, with $N\geq 1$, $I=(0,t^*)$ is the maximal interval of existence  of $u$, $m>1$ and $u_0(\nx)$ is a nonngative and sufficiently regular function.  The problem is equipped with different  boundary conditions and depending on such boundary conditions as well as on the expression of the source $g$, global existence and blow-up  criteria for solutions to \eqref{ProblemAbstract} are established. Additionally, in the three dimensional setting and when blow-up occurs, lower bounds for the blow-up time $t^*$ are also derived.  
\end{abstract}
\section{Introduction and motivations}\label{IntroductionSection} 
It is well known that several natural phenomena appearing in various physical, chemical and biological applications, are modelled through reaction diffusion equations. Their description, generally given in a cylinder $\Omega\times I$, where $\Omega$ is a bounded smooth domain  of $\R^N$ ($N\geq 1$) with regular boundary $\partial \Omega$, and $I=(0,t^*)$, is formulated by an initial boundary value problem in the unknown $u=u({\bf x},t)$ reading as
\begin{equation}\label{ReacDiffsytem}
\begin{cases}
u_t=\nabla \cdot \mathcal{A}(u,\nabla u,{\bf x},t)+\mathcal{B}(u,\nabla u,{\bf x},t)& {\bf x} \in \Omega, t\in I,\\
u({\bf x},0)=u_0({\bf x})& {\bf x} \in  \Omega,\\
\textrm{Boundary conditions on $u$} &  {\bf x} \in \partial \Omega, t\in I.
\end{cases}
\end{equation}
As to the question tied to the existence of local (i.e., $t^*$ finite) or global (i.e., $t^*=\infty$) solutions to classes of nonlinear problems of this type,  sufficient conditions on  $\mathcal{A}$ (as for instance, standard ellipticity behavior) as well as growth and regularity assumptions on both $\mathcal{A}$ and $\mathcal{B}$ guaranteeing this existence are known and have been widely studied in the literature (we refer, for instance,  to \cite{Ball77remarkson,Kielhofer1974,Krylov,LSUBookInequality}). 

In this paper we dedicate our attention to  problem \eqref{ReacDiffsytem} in the case  $\mathcal{A}(u,\nabla u,{\bf x},t)=\nabla u^m$ and $\mathcal{B}(u,\nabla u,{\bf x},t)=g(u,|\nabla u|)$ and endowed with some boundary conditions, i.e.
\begin{equation}\label{General_Problem}
\begin{cases}
u_t=  \Delta u^m+ g(u,|\nabla u|) &  {\bf x}\in \Omega, t\in I, \\
k u_{\boldsymbol\nu}+h u=0 & {\bf x}\in \partial\Omega, t \in  I, \\
u({\bf x},0)=u_0({\bf x}) & {\bf x} \in\Omega, 
\end{cases}
\end{equation}
where $\Omega$ and $I$ were already introduced in the description of \eqref{ReacDiffsytem}. Further, ${\bf \boldsymbol\nu}=(\nu_1,\ldots,\nu_N)$ stands for the outward normal unit vector to the boundary $\partial \Omega$,  $u_{ \boldsymbol\nu}$ is the normal derivative of $u$, $m>1$, $k\geq 0$ and $h>0$. Additionally, $u_0:=u_0({\bf x})\not\equiv 0$ is a nonnegative sufficiently smooth function (possibly also verifying compatible conditions on $\partial \Omega$), and $g(u,|\nabla u|)$ is a regular function of its arguments  and is such that $\underline{u}\equiv 0$ represents a subsolution of the first equation in \eqref{General_Problem}; henceforth, through the maximum principle, the nonnegativity on $\Omega \times I$ of solutions $u$ to  \eqref{General_Problem} remains essentially justified (see \cite{LSUBookInequality,ProtMurrWeinberger}).

Beyond problems arising in the mathematical models for gas or fluid flow in porous media (see \cite{Aronson1986} and \cite{vazquez2007porous}), the formulation in \eqref{General_Problem} also describes  the evolution of some biological population $u$ occupying a certain domain whose growth is governed by the law of $g$ (see \cite{GURTIN197735}); precisely, the term $\Delta u^m$ idealizes the spread of the population, the parameter $m$ indicating the speed of propagation: $m>1$ corresponds to slow, $0<m<1$ fast and the limit case $m=1$ infinity propagation. Moreover, when the coefficient $k$ is zero (the well known Dirichlet boundary conditions), then the distribution of $u$ on the boundary  of the domain maintains constant through the time, while for $k,h>0$ the Robin boundary conditions are recovered: they model a negative flux on the boundary, virtually meaning that the population $u$ gets out of the domain with rate $-h/k$. 

There are several investigations concerning different variants of the initial bo\-und\-ary value problem \eqref{General_Problem}, all devoted to existence and properties of solutions: global and/or local existence, lower and upper bound of blow-up time, blow-up rates and/or asymptotic behavior. In our opinion, the following papers deserve to be referred also because they inspire this present work.
\begin{itemize}
\item \textit{Linear diffusion case ($m=1$) and  $g(u,|\nabla u|)= u^p$, with $p>1$}. For $\Omega=\R^N$, $N\geq 1$, in \cite{ARONSON197833}, \cite{Fujita_1966}  and \cite{kobayashi1977} it is shown that for $1<p\le 1+(2/N)$ the problem has no global positive solution, whilst for $p>1+(2/N)$ it is possible to fix appropriate initial data $u_0$ emanating global solutions. When $\Omega$ is a bounded and smooth domain of $\R^3$ and Dirichlet boundary conditions are assigned, in \cite{PAYNE_Shaefer_2007} a lower bound for the blow-up time of solutions, if blow-up occurs, is derived, and \cite{PayneSchaefer_Robin} essentially deals with blow-up and global existence questions for the same problem in the $N$-dimensional setting, with $N\geq 2$, and endowed with Robin boundary conditions. 
\item \textit{Linear diffusion case ($m=1$) and  $g(u,|\nabla u|)= k_1u^p-k_2|\nabla u|^q$, $k_1,k_2>0$ and $p,q\geq 1$}. In \cite{Souplet_Gradient} it is proved that for  $q = 2p/(p+1)$ and small $k_2>0$ blow-up can occur for any  $N \geq 1$, $p > 1$, $(N - 2)p < N + 2$ and without any restriction on the initial data, while lower bounds of the blow-up time, if blow-up occurs, are derived in \cite{MarrasViglialoroVernier_Kodai} when $k_1$ and $k_2$ are time dependent functions and under different boundary conditions. 
\item \textit{Nonlinear diffusion case ($m>1$) and  $g(u,|\nabla u|)=u^p$, with $p>1$}. For $\Omega=\R^N$,  $N\geq 1$, in  \cite{galaktionov_1994}, \cite{GalaktionovKurdyMikha} and  \cite{LevineTheRoleOf} it is shown that for $1<p\le m+(2/N)$ the problem has no global positive solution, whilst for $p>m+(2/N)$ there exist initial data $u_0$ emanating global solutions. When $\Omega$ is a bounded and smooth domain of $\R^N$, $N\geq 1$, and under Dirichlet boundary conditions, in \cite{Galaktionov1981_Russian} is proved that for $1<p<m$ the problem admits global solutions for all $u_0$ such that $u_0^{m-1}\in H_0^1(\Omega),$ while for $m<p<m(1+(2/N))+(2/N)$ specific initial data produce unbounded solutions (see also \cite{Sacks}).  It is also worth to mention that \cite{GalaktionovShmarevVazquez} and \cite{GalaktionovVazquezExtinction} focus on results dealing with regularity and asymptotic behavior of solutions when $g(u,|\nabla u|)=-u^p$, with $p>0$, defined in the whole space $\R^N$, with $N\geq 1$. 
\item \textit{Nonlinear diffusion case ($m>1$) and  $g(u,|\nabla u|)= u^p-u^\mu |\nabla u^\alpha|^q$, with $p,q,\alpha\geq  1$ and $\mu\geq 0$}. With $\Omega$ bounded and smooth in $\R^N$, $N\geq 1$, and under Dirichlet boundary conditions, in \cite{AndreuEtAl} the authors treat the existence of the so called \textit{admissible solutions} and  show that they are globally bounded  if   $p<\mu + mq$ or $m< p=  \mu + mq$, as well as  the existence of blowing up admissible solutions, under the complementary
condition $1\leq \mu  + mq < p$. Similarly, for $\alpha=m+\mu/q$, $m \geq 1, m/2 + \mu/q > 0, 1  \leq q < 2$, existence of global \textit{weak solutions} is addressed in \cite{AndrMazSimoToledo}.
\end{itemize}
In the context of this premise, we remark that our investigation is not focused on the question concerning the existence of solutions to system \eqref{General_Problem}, but rather on their maximal interval of existence $I$.  In particular, in the framework of nonnegative \textit{classical solutions}, we follow the same approach used in largely cited papers (see, for instance, \cite{PaynPhilSchaefer,PAYNE_Shaefer_2007,PaynePhiProytc,PhilippProyt,ShaeferWithoutGradient,ShaeferExistClassicaPorous} and references therein, for linear or nonlinear diffusion equations, even including our same case, i.e. systems like \eqref{General_Problem} with $m>1$) where such an existence is a priori assumed. Additionally, as to the lifespan $I$ of these solutions, only two scenarios can appear and they provide the following \textit{extensibility criterion} (\cite{Ball77remarkson,BandleBrunnerSurvey,Kielhofer1974}):  
\begin{equation}\label{ExtensibilityCrit}
\begin{aligned}
& \triangleright)\; I=(0,\infty), \textrm{so that $u$ remains bounded for all $\nx\in\Omega$ and time $t>0$}, \\
& \triangleright)\; I=(0,t^*), \textrm{$t^*$ finite (the blow-up time),  so that  $\lim_{t \rightarrow t^*}\lVert u(\cdot,t)\rVert_{L^\infty(\Omega)} =\infty$.} 
\end{aligned}
\end{equation}
By analyzing the expressions of $g$ presented in the previous items, it is reasonable to expect that the contribution of the positive power addendum, representing a source which essentially increases the energy of the system, stimulates the occurrence of the blow-up; conversely,  the negative terms have a damping effect, absorbs the energy and, so,  contrasts the power source term.

Exactly in line with the state of the art above reviewed, with this paper we aim at expanding the underpinning theory of the mathematical analysis for problem \eqref{General_Problem} when different choices of $g$, $h$ and $k$ are considered. Indeed, to the best of our knowledge, the interplay between both positive and negative powers of $u$, or $|\nabla u|$, in the source $g$ and the Robin/Dirichlet boundary conditions has not yet been extensively studied. To be precise, our contribution includes blow-up and global existence criteria for \textit{nonnegative and classical solutions} to \eqref{General_Problem} and estimates of the blow-up time when it occurs. We proof three theorems that are summarized as follows:
\begin{itemize}
\item [$\bullet$] \textit{Criterion for blow-up in $\R^N$, $N\geq 1$: Theorem \ref{TheoremBlowUpDifferencePower}.}  If $g(u,|\nabla u|)=k_1u^p-k_2u^q$, $k_1,k_2,h>0$, $k=1$, $p\geq\max\{m,q\}$ with $m,q> 1$, and $u_0({\bf x})$ is large enough, then the lifespan $I$ of the nonnegative classical solution of problem \eqref{General_Problem}  is finite and $u$ blows-up at some finite time $t^*$. 
\item [$\bullet$] \textit{Criterion for global existence in $\R^N$, $N\geq 1$: Theorem \ref{TheoremGlobalferencePower}.} If $g(u,|\nabla u|)=k_1u^p-k_2u^q$, $k_1,k_2,h>0$, $k=1$, $p<m$ with $m,q> 1$, then for any $u_0(\nx)$ the lifespan $I$ of the nonnegative classical solution of problem \eqref{General_Problem}  is infinite and $u$ is bounded for all time $t>0$. 
\item [$\bullet$] \textit{Lower bound of the blow-up time in $\R^3$: Theorem \ref{TheoremEstimateBlowUp}}. If $g(u,|\nabla u|)=k_1u^p-k_2|\nabla u|^q$, $k_1,k_2,h>0$, $k=0$, for  $p\geq 2$, $2-1/p<m<p$ and $p\geq q \geq 2$, and $u$ is a nonnegative classical solution of problem \eqref{General_Problem} which becomes unbounded in a certain measure and at some finite time $t^{*}$, then, if $k_2$ is sufficiently large, there exists $T$ such that $t^{*}\geq T$.  
\end{itemize}
\begin{remark}
Even if the main motivation of this paper lies in enhancing the mathematical theory tied to nonlinear partial differential equations,  we want to underline that the expressions of the function $g$ given above are justified also by applicative reasons. Indeed, according to \cite{Souplet_Gradient}, a single (biological) species density $u$ occupying a bounded portion of the space evolves in time by displacement, birth/reproduction and death. In particular, the births are described by a superlinear power of such a distribution, the natural deaths by a linear one and the accidental deaths by a function of its gradient; it leads to $u_t=\Delta u +C_1u^p-C_2u-C_3 |\nabla u|^q$, with $p,q>1$ and $C_1,C_2,C_3 >0$. Adding to this equation homogeneous Dirichlet conditions corresponds to a non-viable environment on the boundary; homogeneous Neumann conditions stand for a totally insulated domain and Robin ones to a domain which allows the species to cross the boundary. Furthermore, other models originally introduced for only a single species describe the population growth through the preceding equation in which the source $C_1u^p-C_2u-C_3 |\nabla u|^q$ is replaced by the so called logistic function $u(a-bu)$, with $a,b>0$ (\cite{Verhulst}, Pierre-François Verhulst (1804--1849)), or more generally by functions independent of the accidental deaths and whose qualitative behavior is $u^l(1-u)$, with $l\ge 1$. All the mentioned sources have been also employed in chemotaxis models, precisely to describe the self-organizing of living organisms (\cite{AiTsEdYaMi,CaoZheng,ViglialoroDifferentialIntegralEquations,ViglialoroBoundnessVeryWeak,ViglialoroWolleyDCDS,WinklerVeryWeakSol,WinDespiteLogistic}).
\end{remark}
\section{Main assumptions and preparatory lemmas}
In this section we give some crucial hypothesis and lemmas which will be considered through the paper in the proofs of the main theorems. 

First we give these
\begin{state}
Let $m>1$, $h>0$, $p,q\geq 1$ be real numbers and $N\geq 1$; we establish that  
\begin{enumerate}[label=($\mathcal{H}$$_\arabic*$)]
\item\label{Assumption0}
 For $2-1/p<m<p$, $p\geq q \geq 2$  
\begin{equation}\label{DefinitionCoeff}
\begin{cases}
s=p-1, \quad \mu=\frac{q-1}{s}<1, \quad  d=\frac{m-1}{s}<1,\\ 
1 < \delta< \frac{2}{3}(m+d)\frac{2m+3d-3}{2m+3d-1},  \quad \alpha= \frac{2(m+d)-\delta}{2(m+d)-3\delta}>1,\\
 \beta=\frac{2m+3d-1}{2m+3d-3\alpha}>1, \quad \sigma=\frac{2(m+d)-3\delta}{2(m+d)}>0, \quad \gamma=d+\delta>1.
\end{cases}
\end{equation}
\item\label{Assumption1} $\Omega$ is a bounded domain of $\mathbb{R}^N$, star-shaped and convex in two orthogonal directions, whose geometry for some origin ${\bf x}_0$ inside $\Omega$ is such that 
\begin{equation}\label{GeometricProperties} 
m_1:= \frac{3}{2 \min_{\partial \Omega} (({\bf x}- {\bf x}_0) \cdot {\bf \boldsymbol\nu} )}>0,\quad m_2:= 1+ \frac{\max_{\overline{\Omega}} | {\bf x}- {\bf x}_0|}{\min_{\partial \Omega} (({\bf x}- {\bf x}_0) \cdot {\bf \boldsymbol\nu})}>1.
\end{equation}
\item\label{Assumption2} $\Omega$ is a bounded smooth domain of $\mathbb{R}^N$ such that 
 \begin{equation}\label{Condition_First_Eigen_Rob}
 \frac{\xi_1(hm)}{hm} \geq \frac{2m_1N}{3}+m_2-1,
 \end{equation} 
being $m_1$ and $m_2$ as in \eqref{GeometricProperties} and $\xi_1(h)$ the first positive eigenvalue associated to the supported membrane problem 
\begin{equation}\label{FreeMembraneProblem}
  \begin{cases} 
 \Delta w + \xi (h) w=0& {\bf x}\,\in \,  \Omega, \\
 w_{ \boldsymbol\nu}+h w=0&{\bf x}\,\in\, \partial \Omega.\\
\end{cases}
 \end{equation}
 \end{enumerate}
\end{state}
The forthcoming two lemmas will be employed in the proof of Theorem \ref{TheoremGlobalferencePower} exactly in order to estimate a certain non-zero boundary integral  when Robin boundary conditions are considered in system \eqref{General_Problem}. In particular, Lemma \ref{Lemma_Tecnico}  uses the result in Lemma \ref{lemma1} and even though it was already derived in \cite[Lemma 3.3]{ViglialoroDifferentialIntegralEquations}, for the sake of completeness  we include its proof. 
\begin{lemma}\label{lemma1} 
Let $\Omega$ be a domain of $\mathbb{R}^N$ verifying assumption \ref{Assumption1}. For any nonnegative $C^1(\bar{\Omega})$-function $V$, we have
\begin{equation}\label{SobolevTypeInequBoundary}
\int_{\partial \Omega} V^2 d s \leq \frac{2m_1 N}{3}\int_\Omega V^2 d {\bf x}+2(m_2-1)\int_\Omega V|\nabla V|d {\bf x}.
 \end{equation}
\begin{proof}
This is relation \cite[(A.1) of Lemma A.1.]{PaynePhilVernierSECOND} with $n=1$ and written in terms of the coefficients in \eqref{GeometricProperties} of \ref{Assumption1}. 
\end{proof}
\end{lemma}
\begin{lemma}\label{Lemma_Tecnico}  
Let $\Omega$ be a domain of $\mathbb{R}^N$ verifying assumptions \ref{Assumption1} and \ref{Assumption2}. For any  nonnegative $C^1(\bar{\Omega})$-function $V$ verifying $V_{ \boldsymbol\nu}+h V=0$ on $\partial \Omega$, we have
 \begin{equation}\label{IneqLemma_for_xi}
 \int_\Omega |\nabla V^m|^{2} d{\bf x}\geq \sigma \int_\Omega V^{2m}d{\bf x},
 \end{equation}
 with $\sigma:=\eta(hm)$, being $\eta(h)=(3\xi_1(h)-h(2m_1N+3m_2-3))/(3(h(m_2-1)+1))$.
\begin{proof}
For any nonnegative $C^1(\bar{\Omega})$-function $V$ such that $V_{ \boldsymbol\nu}+h V=0$ on $\partial \Omega$, the general Poincaré inequality returns this relation for the first eigenvalue $\xi_1(h)$ of \eqref{FreeMembraneProblem}:
\begin{equation*} 
\xi_1(h)\int_\Omega V^2 d{\bf x} \leq \int_\Omega |\nabla V|^2d{\bf x}+h\int_{\partial \Omega}V^2d s.
\end{equation*} 
It can be written, through relation \eqref{SobolevTypeInequBoundary} and subsequently the Young inequality with exponents $1/2$, as
\begin{equation*}
\xi_1(h)\int_{\Omega} V^{2} d{\bf x} \leq h\Big(\frac{2m_1N}{3}+m_2-1\Big)\int_{\Omega} V^{2}d{\bf x}+(h(m_2-1)+1)\int_\Omega |\nabla V|^2 d{\bf x}.
\end{equation*}
Hence, since $m_2>1$, we also have
\begin{equation}\label{IneqLemma_1}
\int_{\Omega} |\nabla V|^2 d{\bf x} \geq \eta (h) \int_\Omega  V^2 d{\bf x},
\end{equation}
with $\eta(h)=(3\xi_1(h)-h(2m_1N+3m_2-3))/(3(h(m_2-1)+1)).$ The function $\varphi=V^m$ verifies $\varphi_{ \boldsymbol\nu}+hm \varphi=0$, so that \eqref{IneqLemma_1} provides the 
 \begin{equation*} 
\int_\Omega |\nabla V^m|^2 d{\bf x}=\int_\Omega |\nabla \varphi|^2 d{\bf x} \geq \eta(hm) \int_\Omega \varphi^2 d{\bf x}=\sigma \int_\Omega V^{2m} d{\bf x},
 \end{equation*} 
 with $\sigma:=\eta(hm)=(\xi_1(hm)-hm(2m_1N/3+m_2-1))/(hm(m_2-1)+1)$, that is nonnegative by  \eqref{Condition_First_Eigen_Rob}.
  \end{proof}
  \end{lemma}
  Similarly, the remaining two lemmas of this section are necessary to arrange terms emerging in the proof of Theorem \ref{TheoremEstimateBlowUp}, some of them also depending on $|\nabla u|$.
  \begin{lemma}\label{LemmaInR3Embedding}
Let $\Omega$ be a bounded domain of $\mathbb{R}^3$ with Lipschitz boundary, and let  $m, \alpha, \beta, d, \sigma$ and $\gamma$ the constants defined in \eqref{DefinitionCoeff} of \ref{Assumption0}.  Then for any  nonnegative $C^1(\bar{\Omega})$-function $V$
\begin{equation}\label{BoundVm+1}
 \int_{\Omega} V^{m+1}d{\bf x}\leq  \frac{\gamma - 1} {\gamma -  \mu}\varepsilon_1  \int_{\Omega}  V^{m+\mu} d{\bf x} +\frac{1 - \mu} {\gamma -  \mu}  \varepsilon_1^{-\frac{\gamma - 1} {1 -  \mu}}\int_{\Omega} V^{m+\gamma}d{\bf x},
\end{equation} 
where $\varepsilon_1$ is an arbitrary positive constant.

If, additionally, $V$ vanishes on $\partial \Omega$, there exists a positive $\Gamma$ such that for every $\varepsilon_2>0$
\begin{equation}\label{BoundVm+DetlaFinal}  
\begin{split}
\int_\Omega V^{m+\gamma}d{\bf x}\leq  & \Gamma^{\frac{3\delta}{m+d}+\frac{6\alpha}{2m+3d}}3\alpha d\frac{\sigma}{2m+3d}
\varepsilon_2\int_\Omega |\nabla V^{\frac{m+d}{2}}|^2 d{\bf x} \\&
+\Gamma^{\frac{3\delta}{m+d}+\frac{6\alpha}{2m+3d}}d\sigma \frac{2m+3d-3\alpha}{2m+3d}
\varepsilon_2\Big(\int_\Omega V^m d{\bf x}\Big)^{\alpha\beta}
\\&
+(1-d)\varepsilon_2 \Gamma^\frac{3\delta}{m+d}\sigma
\Big(\int_\Omega V^m d{\bf x}\Big)^{\alpha}\\&+\Gamma^\frac{3\delta}{m+d}\frac{3\delta}{2(m+d)}\varepsilon_2^{1-\frac{2(m+d)}{3\delta}}\int_\Omega |\nabla V^{\frac{m+d}{2}}|^2 d{\bf x}.
\end{split}
\end{equation}
\begin{proof}
Since $\mu<1$, for some positive constant $\gamma>1$ the Young inequality and the consideration of $\varepsilon_1>0$ yield 
\begin{equation*}
\begin{split}
\int_{\Omega} V^{m + 1}d{\bf x} &\leq \Big( \int_{\Omega} V^{m + \mu}d{\bf x}\Big)^{\frac{\gamma - 1} {\gamma -  \mu}}\Big( \int_{\Omega} V^{m+\gamma }d{\bf x}\Big)^{\frac{1 - \mu} {\gamma -  \mu}} \\ &
\quad  \frac{\gamma - 1} {\gamma -  \mu}\varepsilon_1  \int_{\Omega}  V^{m+\mu} d{\bf x} +\frac{1 - \mu} {\gamma -  \mu}  \varepsilon_1^{-\frac{\gamma - 1} {1 -  \mu}}\int_{\Omega} V^{m+\gamma}d{\bf x},
\end{split}
\end{equation*} 
so that the first thesis is shown.

On the other hand, let $V$ be such that $V=0$ on $\partial \Omega$: the Sobolev embedding in $\R^3$, $W_0^{1,2} \hookrightarrow L^6$, 
 provides
 \begin{equation}\label{SobolevInequ3D}
 \int_\Omega \Big(V^\frac{m+d}{2}\Big)^6d{\bf x}\leq \Gamma^6 \Big(\int_\Omega |\nabla V^\frac{m+d}{2}|^2d{\bf x}\Big)^3,
 \end{equation}
$\Gamma = 4^{1/2} 3^{-1/2}\pi^{-2/3}$ being the best Sobolev constant (see \cite{Talenti1976BestConstant}). Now, for $\gamma=d+\delta>1$, the H\"{o}lder inequality leads to 
 \begin{equation}\label{HolderConDelta}
 \begin{split}
 \int_{\Omega} V^{m + \gamma}d{\bf x} & =\int_{\Omega} V^{(m + d)+\delta}d{\bf x} \\ & \leq \Big( \int_{\Omega} V^{m + d}d{\bf x}\Big)^{\frac{2(m+d)-\delta} {2(m+d)}}\Big( \int_{\Omega} (V^\frac{m+d}{2})^6d{\bf x}\Big)^{\frac{\delta} {2(m+d)}},
 \end{split}
 \end{equation} 
so that by replacing \eqref{SobolevInequ3D} into \eqref{HolderConDelta}, we obtain
 \begin{equation*}
  \int_{\Omega} V^{(m + d)+\delta}d{\bf x} \leq \Gamma^{\frac{3\delta}{m+d}}\Big( \int_{\Omega} V^{m + d}d{\bf x}\Big)^{\frac{2(m+d)-\delta} {2(m+d)}}\Big(\int_{\Omega} |\nabla V^\frac{m+d}{2}|^2d{\bf x}\Big)^{\frac{3\delta} {2(m+d)}}.
  \end{equation*}
The introduction of an arbitrary positive constant $\varepsilon_2$, and an application of the Young inequality, allow us to write   (recall \eqref{DefinitionCoeff})
 \begin{equation}\label{BoundVm+Detla}
 \begin{split}
   \int_{\Omega} V^{(m + d)+\delta}d{\bf x}& \leq 
   \Gamma^{\frac{3\delta}{m+d}}\Big(\varepsilon_2\Big( \int_{\Omega} V^{m + d}d{\bf x}\Big)^{\frac{2(m+d)-\delta} {2(m+d)-3\delta}}\Big)^{\frac{2(m+d)-3\delta} {2(m+d)}} \\ & 
   \quad \times \Big(\varepsilon_2^{1-\frac{2(m+d)}{3\delta}} \int_{\Omega} |\nabla V^\frac{m+d}{2}|^2d{\bf x}\Big)^\frac{3\delta} {2(m+d)}
   \\&
   \leq \Gamma^{\frac{3\delta}{m+d}}\varepsilon_2{\frac{2(m+d)-3\delta} {2(m+d)}}\Big( \int_{\Omega} V^{m + d}d{\bf x}\Big)^{\frac{2(m+d)-\delta} {2(m+d)-3\delta}}\\ &  \quad +\Gamma^{\frac{3\delta}{m+d}}\varepsilon_2^{1-\frac{2(m+d)}{3\delta}}\frac{3\delta} {2(m+d)} \int_{\Omega} |\nabla V^\frac{m+d}{2}|^2d{\bf x}.
   \end{split}
   \end{equation}
To bound the term $\Big(\int_{\Omega} V^{m + d}d{\bf x}\Big)^{(2(m+d)-\delta)/(2(m+d)-3\delta)}$, let us observe that the H\"{o}lder and the Schwarz inequalities give, respectively, 
\begin{equation}\label{BoundVM+1}
\int_\Omega V^{m+1}d{\bf x} \leq \Big(\int_\Omega V^{2(m+d)}d{\bf x}\Big)^\frac{1}{m+2d}\Big(\int_\Omega V^m d{\bf x}\Big)^\frac{m+2d-1}{m+2d},
\end{equation}
and 
\begin{equation}\label{BoundV2M+D}
\int_\Omega V^{2(m+d)}d{\bf x} \leq \Big[\int_\Omega \Big(V^\frac{m+d}{2}\Big)^6d{\bf x}\int_\Omega V^{m+d}d{\bf x}\Big]^\frac{1}{2}.
\end{equation}
Now, using in \eqref{BoundV2M+D} relation \eqref{SobolevInequ3D}, we get
\begin{equation*}
\int_\Omega V^{2(m+d)}d{\bf x} \leq \Gamma^3 \Big(\int_\Omega |\nabla V^{\frac{m+d}{2}}|^2d{\bf x}\Big)^\frac{3}{2}\Big(\int_\Omega V^{m+d}d{\bf x}\Big)^\frac{1}{2},
\end{equation*}
and hence \eqref{BoundVM+1} reads
\begin{equation}\label{A}
\begin{split}
 \int_\Omega V^{m+1}d{\bf x} & \leq  \Gamma^\frac{3}{m+2d}\Big(\int_\Omega |\nabla V^{\frac{m+d}{2}}|^2d{\bf x}\Big)^\frac{3}{2(m+2d)}  \Big(\int_\Omega V^{m+d}d{\bf x}\Big)^{\frac{1}{2(m+2d)}} \\&
\quad \times \Big(\int_\Omega V^m d{\bf x}\Big)^{\frac{m+2d-1}{m+2d}}.
\end{split}
\end{equation}
In addition, we first use again the H\"{o}lder inequality to lead to
\begin{equation}\label{AA}
\int_\Omega V^{m+d}d{\bf x} \leq \Big(\int_\Omega V^{m+1}d{\bf x}\Big)^d \Big(\int_\Omega V^m d{\bf x}\Big)^{1-d},
\end{equation}
and then we insert this estimate in \eqref{A}; combining terms, applying 
\begin{equation}\label{(eq:2.20)}
  a^r b^{1-r}\leq r a + (1-r) b,
\end{equation} 
valid for $a,b\geq 0$ and $ 0<r<1$,
we arrive at ($\alpha$ as in \eqref{DefinitionCoeff} of \ref{Assumption0})
\begin{equation}\label{EqVEm+d}
\begin{split}
\Big(\int_\Omega V^{m+1}d{\bf x}\Big)^\alpha & \leq \Gamma^{\frac{6\alpha}{2m+3d}}\Big(\int_\Omega |\nabla V^\frac{m+d}{2}|^2d{\bf x}\Big)^\frac{3\alpha}{2m+3d}\Big(\int_\Omega V^m d{\bf x}\Big)^{\alpha\frac{2m+3d-1}{2m+3d}}\\&
\leq \Gamma^{\frac{6\alpha}{2m+3d}} \frac{3\alpha}{2m+3d} \int_\Omega |\nabla V^\frac{m+d}{2}|^2d{\bf x}\\ & 
\quad + \Gamma^{\frac{6\alpha}{2m+3d}} \frac{2m+3d-3\alpha}{2m+3d}\Big(\int_\Omega V^m d{\bf x}\Big)^{\alpha\frac{2m+3d-1}{2m+3d-3\alpha}}.
\end{split}
\end{equation}
Hence, by rearranging again \eqref{AA} with \eqref{(eq:2.20)} we attain
\begin{equation*}
\begin{split}
\Big(\int_\Omega V^{m+d}d{\bf x}\Big)^\alpha & \leq \Big[\Big(\int_\Omega V^{m+1} d{\bf x} \Big)^d \Big(\int_\Omega V^m d{\bf x}\Big)^{1-d}\Big]^\alpha
\\ & \leq d \Big(\int_\Omega V^{m+1} d{\bf x} \Big)^\alpha+(1-d)\Big(\int_\Omega V^m d{\bf x}\Big)^\alpha,
\end{split}
\end{equation*}
so that in view of \eqref{EqVEm+d} expression \eqref{BoundVm+Detla}  (recall $\gamma=d+\delta$) infers our thesis.
\end{proof}
  \end{lemma}
  \begin{lemma}\label{PropTec} 
  Let $m,d$ and $\delta$ as in \eqref{DefinitionCoeff} of  \ref{Assumption0}. If $c_1, c_2,\ldots,c_6$ are positive real numbers satisfying
  \begin{equation}\label{ConditionCiProp}
  c_3 \geq c_5\Big(\frac{c_5}{c_6}\Big)^\frac{-3\delta}{2(m+d)}\Big(\frac{3\delta}{(2m+2d-3\delta)}\Big)^{-\frac{3\delta}{2(m+d)}}\frac{2(m+d)}{2(m+d)-3\delta},
  \end{equation}
  then there exits $\xi_m \in (0,\infty)$ such that 
  \begin{equation}\label{PropInequality}
  c_5\xi_m+c_6\xi_{m}^{1-\frac{2(m+d)}{3\delta}}-c_3\leq 0.
  \end{equation}
  \begin{proof}
  For any $\xi \in (0,\infty)$, the function $\Phi(\xi):=c_5\xi+c_6\xi^{1-2(m+d)/3\delta}$ attains its minimum at the point 
  \begin{equation}\label{MinimumPointFunctionPhi}
 \xi_m=\Big(\frac{3\delta c_5}{c_6(2m+2d-3\delta)}\Big)^\frac{-3\delta}{2(m+d)}.
  \end{equation}
  Therefore, since \eqref{ConditionCiProp} holds we have
  \[ c_3\geq c_5\Big(\frac{c_5}{c_6}\Big)^\frac{-3\delta}{2(m+d)}\Big(\frac{3\delta}{(2m+2d-3\delta)}\Big)^{-\frac{3\delta}{2(m+d)}}\frac{2(m+d)}{2(m+d)-3\delta}=\Phi(\xi_m),\]
  and relation \eqref{PropInequality} is proven.
  \end{proof}
    \end{lemma}
\section{Analysis and proofs of the main results}
In this section we discuss and give the demonstrations of our main theorems, whose general overview was summarized in $\S$\ref{IntroductionSection}. 
  \subsection{A criterion for blow-up}
The first theorem is dedicated to understand properties of solutions to system  \eqref{General_Problem} when $g(u,|\nabla u|)=k_1u^p-k_2u^q$ and under Robin boundary conditions. Essentially, we observe that if the power $q$ of the absorption term in $g$, as well as the coefficient $m$ of the diffusion, do not surpass the power $p$ from the growth contribution, then the occurrence of blow-up phenomena at some finite time may appear for some initial data $u_0({\bf x})$, despite the outflow boundary conditions; in particular no global solution is expected. 
\begin{theorem}\label{TheoremBlowUpDifferencePower}
Let $\Omega$ be a bounded smooth domain of $\R^N$, $N\geq 1,$  $k_1,k_2,h>0$, $k=1$, $p\geq\max\{m,q\}$, with $m,q > 1$, $g(u,|\nabla u|)=k_1u^p-k_2u^q$ and $u_0({\bf x})\not \equiv 0$ a nonnegative function from $C^1(\bar{\Omega})$. Moreover, let  $u\in C^{2,1}(\Omega \times (0,t^*))\cap C^{1,0}(\bar{\Omega} \times [0,t^*))$ be the nonnegative solution of problem \eqref{General_Problem}. If 
\begin{equation*} 
\begin{split}
\psi(t):&=-\frac{p+m}{2m}\int_\Omega \lvert \nabla u^m\rvert^2 d \nx + k_1 \io u^{p+m} d \nx  \\ &
\quad  -k_2\io u^{q+m} d \nx-\frac{h(p+m)}{2}\int_{\partial \Omega} u^{2m}ds, \quad \textrm{for all $t$}  \in (0,t^*),
\end{split}
\end{equation*}
is such that $\psi(0)> 0$, then $t^*<\infty$, or equivalently $I=(0,t^*)$. In particular,  $\lVert u(\cdot,t)\rVert_{L^\infty(\Omega)}\nearrow \infty$ as $t\searrow t^*$ 
at some time $t^*$ satisfying
\begin{equation*}
t^* <T= \frac{m+1}{p-1}\frac{\varphi(0)}{\psi(0)}.
\end{equation*}
\begin{proof}
Let $u$ be the nonnegative classical solution of \eqref{General_Problem} satisfying $u_{ \boldsymbol\nu}=-h u$ on $\partial \Omega$. By a differentiation we can write
\begin{equation}\label{evolutionPsi}
\begin{split}
\frac{1}{m+1}\psi'(t)&=-\frac{p+m}{m}\io \nabla u^m\cdot (\nabla u^m)_t d\nx+k_1(p+m)\io u^{p+m-1}u_t d \nx \\ & 
\quad -k_2(q+m)\io u^{q+m-1}u_t d \nx -h(p+m)m\iob u^{2m-1}u_tds\\ &
=-\frac{p+m}{m}\iob (u^m)_t \nabla u^m \cdot \boldsymbol\nu  ds+ \frac{p+m}{m}\io (u^m)_t \Delta u^m d \nx\\ &
\quad +k_1(p+m)\io u^{p+m-1}u_t d \nx  -k_2(q+m)\io u^{q+m-1}u_t d \nx \\ &
\quad -h(p+m)m\iob u^{2m-1}u_tds \\ &
\geq(p+m)\io u^{m-1}u_t(\Delta u^m+k_1u^p-k_2u^q) d \nx \\ &
=(p+m)\io u^{m-1}(u_t)^2 d \nx \geq 0\quad \textrm{for all $t$}\in (0,t^*),
\end{split}
\end{equation}
where we have used the integration by parts formula and the assumption $p\geq q$.

Similarly, as to the evolution of $\varphi(t):=\io u^{m+1}d \nx$, we derive
\begin{equation}\label{evolutionPhi}
\begin{split}
\frac{1}{m+1}\varphi'(t)&=\io  u^m (\Delta u^m+k_1u^p-k_2u^q)d\nx \\ & 
=-\io \nabla u^m\cdot \nabla u^m d \nx  +k_1\io u^{p+m} d \nx-k_2\io u^{q+m} d \nx \\ & 
\quad +\iob u^{m}\nabla u^m \cdot \boldsymbol\nu  ds\\&
=-\io \lvert \nabla u^m\rvert^2 d \nx  +k_1\io u^{p+m} d \nx-k_2\io u^{q+m} d \nx \\ &
\quad -mh\iob u^{2m}ds\\&
\geq-\frac{p+m}{2m}\io \lvert \nabla u^m\rvert^2 d \nx+k_1\io u^{p+m} d \nx-k_2\io u^{q+m} d \nx \\ &
\quad  - \frac{h(p+m)}{2}\iob u^{2m}ds\geq \frac{1}{m+1}\psi(t) \quad \textrm{for all $t$}\in (0,t^*),
\end{split}
\end{equation}
where in this case we relied on the fact that $p\geq m$. Now, the hypothesis $\psi(0)> 0$,  \eqref{evolutionPsi} and \eqref{evolutionPhi} yield 
\begin{equation*}
\psi(t)> 0 \quad \textrm{and}\quad \varphi'(t)>0\quad \textrm{on} \quad (0,t^*). 
\end{equation*}
Since by the Young inequality we have that for all $t\in(0,t^*)$
\begin{equation*}
\frac{1}{m+1}\varphi'(t)=\io u^\frac{m+1}{2} u^\frac{m-1}{2}u_td \nx \leq \bigg(\io u^{m+1}d \nx\bigg)^\frac{1}{2} \bigg(\io u^{m-1}(u_t)^2 \nx\bigg)^\frac{1}{2}, 
\end{equation*}
this implies by virtue of the definition of $\varphi$, in conjunction with \eqref{evolutionPsi} and \eqref{evolutionPhi}, 
\begin{equation*}
\varphi(t)\psi'(t)\geq \frac{m+p}{m+1}\varphi'(t)^2\geq  \frac{m+p}{m+1}\psi(t)\varphi'(t) \quad \textrm{on} \quad (0,t^*),
\end{equation*}
or equivalently
\begin{equation*}
\frac{d}{dt}\big(\psi \varphi^{-\frac{m+p}{m+1}}\big)\geq 0  \quad \textrm{on} \quad (0,t^*).
\end{equation*}
Subsequently, an integration on $(0,t)$ with $t<t^*$ infers, being $\varphi(0)>0$,
\begin{equation*}
\psi(t)\geq \psi(0) \varphi(0)^{-\frac{m+p}{m+1}}\varphi(t)^{\frac{m+p}{m+1}}\quad \textrm{on} \quad (0,t).
\end{equation*}
Finally, recalling \eqref{evolutionPhi}, we have
\begin{equation*}
\varphi'(t)\varphi(t)^{-\frac{m+p}{m+1}}\geq \psi(0) \varphi(0)^{-\frac{m+p}{m+1}}\quad \textrm{on} \quad (0,t),
\end{equation*}
and with  $(m+p)/(m+1)>1$ a further integration leads to 
\begin{equation*}
\frac{1}{\varphi(t)^\frac{p-1}{m+1} }\leq \frac{1}{\varphi(0)^\frac{p-1}{m+1}}-\Big(\frac{p-1}{m+1}\Big)\frac{\psi(0)}{\varphi(0)^{\frac{m+p}{m+1}}}t,
\end{equation*}
that, by virtue of the positivity of $\varphi$ cannot hold for $t\geq T=(m+1)\varphi(0)/(p-1)\psi(0)$. In conclusion, the extensibility criterion \eqref{ExtensibilityCrit} implies that $I=(0,t^*)$, for some $t^*<T$.
\end{proof}
\end{theorem}
  \subsection{A criterion for global existence}
In the next result, we are interested to examine the opposite situation described in Theorem \ref{TheoremBlowUpDifferencePower}. Precisely, by considering in system  \eqref{General_Problem} again $g(u,|\nabla u|)=k_1u^p-k_2u^q$ and Robin boundary conditions, we establish that when the effect of the source (coefficient $p$) is weaker than that of the diffusion (coefficient $m$), the negative flux on the boundary prevents blow-up, even for arbitrary large initial data $u_0({\bf x})$ and any small absorption effect (coefficient $q$).  
  \begin{theorem}\label{TheoremGlobalferencePower}
Let $\Omega$ be a bounded smooth domain of $\R^N$, $N\geq 1,$ satisfying assumptions \ref{Assumption1} and \ref{Assumption2}. Moreover let be $k_1,k_2,h>0$, $k=1$, $q \geq 1$, $p<$m, with $m > 1$, $g(u,|\nabla u|)=k_1u^p-k_2u^q$ and $u_0({\bf x})\not \equiv 0$ a nonnegative function from $C^0(\bar{\Omega})$. Then the nonnegative solution  $u\in C^{2,1}(\Omega \times (0,t^*))\cap C^{1,0}(\bar{\Omega} \times [0,t^*))$ of problem \eqref{General_Problem} is global, or equivalently $I=(0,\infty)$.
\begin{proof}
Let $u$ be the nonnegative classical solution of \eqref{General_Problem} satisfying $u_{ \boldsymbol\nu}=-h u$ on $\partial \Omega$. By differentiating $\varphi(t):=\io u^{m+1} d \nx$ we derive
\begin{equation*}
\begin{split}
\frac{1}{m+1}\varphi'(t)&=\io  u^m (\Delta u^m+k_1u^p-k_2u^q)d\nx \\ & 
=-\io \nabla u^m\cdot \nabla u^m d \nx  +k_1\io u^{p+m} d \nx-k_2\io u^{q+m} d \nx \\ & 
\quad +\iob u^{m}\nabla u^m \cdot \boldsymbol\nu ds\\&
\leq -\io \lvert \nabla u^m\rvert^2 d \nx  +k_1\io u^{p+m} d \nx  \quad \textrm{on}\;(0,t^*), 
\end{split}
\end{equation*}
where we have neglected the last two nonpositive integrals.

On the other hand, since $p<m$, we have thanks to the Young inequality and for some $\varepsilon>0$
\begin{equation*}
k_1\io u^{p+m} d \nx \leq \varepsilon \io u^{2m} d \nx+ C(\varepsilon)|\Omega| \quad \textrm{on}\;(0,t^*),
\end{equation*}
with $C(\varepsilon)=(2m \varepsilon /(p+m)k_1)^{(m+p)/(p-m)}(m-p)/2m>0$ (recall $m>p$). 
Subsequently,  
\begin{equation*}
\begin{split}
\frac{1}{m+1}\varphi'(t)&\leq  -\io \lvert \nabla u^m\rvert^2 d \nx   +\varepsilon \io u^{2m} d \nx+ C(\varepsilon)|\Omega|  \\ & 
\leq -\sigma\io   u^{2m} d \nx +\varepsilon \io u^{2m} d \nx+ C(\varepsilon)|\Omega|\\ & 
=-(\sigma-\varepsilon) \io   u^{2m} d \nx + C(\varepsilon)|\Omega|  \quad \textrm{on}\;(0,t^*), 
\end{split}
\end{equation*}
where we have estimated the integral depending on $\lvert \nabla u^m\rvert^2$ by means of \eqref{IneqLemma_for_xi} of Lemma \ref{Lemma_Tecnico} with, of course, $V=u$. By choosing $\varepsilon=\frac{\sigma}{2}$, and by taking in consideration that an application of the Young inequality infers
\begin{equation*}
-\io u^{2m} d\nx \leq - |\Omega|^\frac{1-m}{1+m} \varphi^\frac{2m}{m+1}\quad \textrm{on}\;(0,t^*),
\end{equation*}
the previous estimate reads
\begin{equation*}
\varphi'(t)\leq -C_0  \varphi^\frac{2m}{m+1}(t)+C_1\quad \textrm{on}\;(0,t^*),
\end{equation*}
where $C_0=(m+1)\sigma  |\Omega|^\frac{1-m}{1+m}/2$ and $C_1=(m+1)C(\varepsilon)|\Omega|$; consequently, ODE comparison arguments justify that
\begin{equation*}
\varphi(t)\leq C:=\max\bigg\{\varphi(0),\bigg(\frac{C_1}{C_0}\bigg)^\frac{m+1}{2m}\bigg\}\quad \textrm{on}\;(0,t^*).
\end{equation*}
Finally, well know extension results for ODE's  with locally Lipschitz continuous right side (see, for instance, \cite{grant2014theoryODE}),  show that $t^*=\infty$; indeed, if $t^*$ were finite, $\varphi (t) \nearrow +\infty$ as $t\searrow t^*$ and it would contradict $\varphi(t)\leq C$ on $ (0,t^*).$ In conclusion, again the extensibility criterion \eqref{ExtensibilityCrit} implies $I=(0,\infty)$.
\end{proof}
\end{theorem}
\begin{remark}
Conversely to the demonstration of Theorem  \ref{TheoremBlowUpDifferencePower}, evidently the proof of this last theorem remains valid also for $k_2=0$, that is in complete absence of absorption terms in $g$. In any case, we preferred to consider the expression of the function $g$ in Theorem  \ref{TheoremGlobalferencePower} as that in Theorem \ref{TheoremBlowUpDifferencePower} exactly to better highlight the different behavior of the corresponding solutions to problem \eqref{General_Problem} despite the same source.
\end{remark}
 \subsection{Lower bounds of the blow-up time}
 This last theorem is concerned with lower bounds of the blow-up time $t^*$ for unbounded solutions to \eqref{General_Problem}, when gradient nonlinearities with absorption effects appear in $g$. More precisely, we define  $g(u,|\nabla u|)=k_1u^p-k_2|\nabla u|^q$ and endow the problem with Dirichlet boundary conditions. We are not aware of general results which straightforwardly infer the existence of unbounded solutions to system \eqref{General_Problem} under these hypothesis; nevertheless, in the spirit of the result derived in Theorem \ref{TheoremBlowUpDifferencePower}, for which blow-up occurs for large initial data and despite negative flux on the boundary,  we understand that also in these circumstances seems reasonable to assume the existence of such blowing-up solutions. 
\begin{theorem}\label{TheoremEstimateBlowUp}
Let $\Omega$ be a bounded domain of $\mathbb{R}^3$ with Lipschitz boundary. Moreover let $k_1,h>0$, $k=0$, $p,q,\alpha,\beta,\mu$ and $\gamma$ as in \eqref{DefinitionCoeff} of \ref{Assumption0}, $g(u,|\nabla u|)=k_1u^p-k_2|\nabla u|^q$ and $u_0({\bf x})\not \equiv 0$ a nonnegative function from $C^0(\bar{\Omega})$, satisfying the compatibility condition $u_0({\bf x})=0$ on $\partial \Omega$. Hence, it is possible to find a positive number $\varSigma$ with the following property: If $k_2$ is a positive real satisfying 
 \begin{equation}\label{AssumptionCoeffieintTheoremBlow-Up}
 k_2\geq k_1(k_1\varSigma)^\frac{1-\mu}{\gamma-1},
 \end{equation}
 and $u\in C^{2,1}(\Omega \times (0,t^*))\cap C^{0}(\bar{\Omega} \times [0,t^*))$ is a nonnegative solution of \eqref{General_Problem} such that $W(t)\nearrow +\infty$ as $t\searrow t^{*}$, with some finite $t^{*}$
and 
 \begin{equation}\label{Auxiliar_function}
  W(t) =\int_{\Omega}u^{m(p-1)} d {\bf x},
  \end{equation} 
 then
 \begin{equation*} 
t^{*}    \geq \frac{W(0)^{-\alpha\beta+1}}{(\mathcal{M}W(0)^{(1-\beta)\alpha}+\mathcal{N})(-\alpha\beta+1)},
 \end{equation*}
  $\mathcal{M}$ and $\mathcal{N}$ being two positive computable constants.
\begin{proof}
Let $u$ be the nonnegative classical solution of \eqref{General_Problem} satisfying $u=0$ on $\partial \Omega$ and $t^*$ be the instant of time where the $W$-measure  \eqref{Auxiliar_function} associated to $u$ becomes unbounded. For $s=p-1$, let us differentiate respect to the time $t$ such $W$-measure. Due to the divergence theorem and the boundary conditions, we obtain
\begin{equation}\label{Derivative_Auxiliar_Function}
\begin{split}
W'(t)& = ms \int_{\Omega} u^{ms-1}[\Delta (u^m)+k_1u^p -k_2|\nabla u|^q] d{\bf x} 
\\ & =-ms \int_{\Omega} \nabla u^{ms-1} \cdot 
\nabla (u^m) d x\\&\;\;\;\;+ msk_1 \int_\Omega u^{s(m+1)}d{\bf x} -msk_2\int_{\Omega} u^{ms-1}|\nabla u|^q d{\bf x}\\&
=-m^2 s(ms-1)\int_\Omega u^{ms-3+m}|\nabla u|^2d{\bf x} \\&\;\;\;\;+msk_1 \int_\Omega u^{s(m+1)}d{\bf x} -msk_2\int_{\Omega} u^{ms-1}|\nabla u|^q d{\bf x} \quad \textrm{on} \quad (0,t^*).
\end{split}
\end{equation}
Now, the assumptions given in  \eqref{DefinitionCoeff} of \ref{Assumption0} imply $ms+q-1>2$, so we can invoke inequality \cite[(2.10)]{PayneElAtl2008} achieving
\begin{equation}\label{(eq:2.13)}
\begin{split}
 msk_2\int_{\Omega}u^{ms -1}|\nabla u|^q d{\bf x} &= msk_2\Big(\frac{q}{ms+q-1}\Big)^q\int_\Omega |\nabla 
 u^{\frac{ms+q-1}{q}}|^q d{\bf x} \\&
 \geq msk_2\Big(\frac{2\sqrt{\lambda_1}}{ms + q-1}\Big)^q \int_{\Omega}u^{ms + q-1}  d{\bf x}\quad \textrm{on} \quad (0,t^*),
 \end{split}
\end{equation} 
where $\lambda_1$ is the optimal Poincaré constant.  

From now on, for simplicity we indicate $u^s=V$ so to have 
\begin{equation}\label{eq:2.11}
 |\nabla V|^2= s^2 u^{2(s-1)}|\nabla u|^2.
\end{equation} 
As a consequence, again due to the positions made in \eqref{DefinitionCoeff} of \ref{Assumption0}, it holds that  $(m-2)+d>0$, so that using $\eqref{(eq:2.13)}$ and $\eqref{eq:2.11}$, relation \eqref{Derivative_Auxiliar_Function} becomes
\begin{equation}\label{(eq:2.14)}
\begin{split}
W'(t)&\leq  -c_1 \int_{\Omega}V^{(m-2)+d}|\nabla  V|^2 d{\bf x}+c_2 \int_{\Omega} V^{m + 1}d{\bf x}\\ &
\quad  -msk_2\Big(\frac{2\sqrt{\lambda_1}}{ms + q-1}\Big)^q  \int_{\Omega} V^{m + \mu}d{\bf x}\\&
=-c_3 \int_{\Omega}|\nabla  V^\frac{m+d}{2}|^2 d{\bf x}+c_2 \int_{\Omega} V^{m + 1}d{\bf x} \\ &
\quad -msk_2\Big(\frac{2\sqrt{\lambda_1}}{ms + q-1}\Big)^q  \int_{\Omega} V^{m + \mu}d{\bf x}\quad \textrm{on} \quad (0,t^*),
 \end{split}
\end{equation} 
where 
\begin{equation*} 
c_1=\frac{m^2(ms-1)}{s},\quad c_2=msk_1,\quad c_3=\frac{4}{(m+d)^2}c_1.
\end{equation*}
Now we are in the position to apply Lemma \ref{LemmaInR3Embedding}: by  using relation \eqref{BoundVm+1}  with $\varepsilon_1=k_2\big(2\sqrt{\lambda_1}/(ms + q-1)\big)^q(\gamma-\mu)/(k_1(\gamma-1))$, \eqref{(eq:2.14)} is simplified to 
\begin{equation*} 
W'(t)\leq -c_3 \int_{\Omega}|\nabla  V^\frac{m+d}{2}|^2 d{\bf x}+c_4 \int_{\Omega} V^{m + \gamma}d{\bf x}\quad \textrm{on} \quad (0,t^*),
\end{equation*}
where 
\begin{equation*} 
c_4=c_2\frac{1 - \mu} {\gamma -  \mu}  \varepsilon_1^{-\frac{\gamma - 1} {1 -  \mu}},
\end{equation*}
whilst rearranging the term $\int_{\Omega} V^{m + \gamma}d{\bf x}$ through \eqref{BoundVm+DetlaFinal} we obtain
\begin{equation}\label{IneqPhiCasiFinal}
\begin{split}
W'(t)&\leq 
\Big(c_5\varepsilon_2+c_6\varepsilon_2^{1-\frac{2(m+d)}{3\delta}}-c_3\Big)
\int_\Omega |\nabla V^{\frac{m+d}{2}}|^2 d{\bf x}\\ &
\quad +\mathcal{M}\Big(\int_\Omega V^m d{\bf x}\Big)^{\alpha}+\mathcal{N}\Big(\int_\Omega V^m d{\bf x}\Big)^{\alpha\beta}\quad \textrm{on} \quad (0,t^*),
\end{split}
\end{equation}
with 
\begin{equation*} 
\begin{cases}
c_5=\Gamma^{\frac{3\delta}{m+d}+\frac{6\alpha}{2m+3d}}\frac{3\alpha dc_4\sigma}{2m+3d},&
c_6= \frac{3\delta c_4 \Gamma^\frac{3\delta}{m+d}}{2(m+d)},\\
\mathcal{M}= \Gamma^\frac{3\delta}{m+d}(1-d)\varepsilon_2\sigma c_4,&
\mathcal{N}= \Gamma^{\frac{3\delta}{m+d}+\frac{6\alpha}{2m+3d}} \frac{2m+3d-3\alpha}{2m+3d}
\varepsilon_2 c_4d\sigma.
\end{cases}
\end{equation*}
Hereafter, setting 
\[\varSigma = 
\frac{\frac{\Gamma^{\frac{3\delta}{m+d}} s^2(m+d)^2}{4 m(ms-1)}}{ \Big[\frac{\gamma-\mu}{\gamma-1}\Big(\frac{2\sqrt{\lambda_1}}{ms+q-1}\Big)^q\Big]^{\frac{\gamma-1}{1-\mu}}}
\frac{1-\mu}{\gamma-\mu}\Big(\frac{6(m+d)\Gamma^{\frac{6\alpha}{2m+3d}}d \alpha \sigma  }{(2m+3d)(2m+2d-3\delta)}\Big)^{1-\frac{3\delta}{2(m+d)}},\]
we observe that using the values of the constants $c_1,c_2,\ldots,c_6$ defined so far, relation \eqref{ConditionCiProp} is precisely equivalent to \eqref{AssumptionCoeffieintTheoremBlow-Up}. Subsequently, Lemma \ref{PropTec} warrants that for $\varepsilon_2=\xi_m$, whose value was computed in \eqref{MinimumPointFunctionPhi}, 
 %
 %
 %
 %
 $c_5\varepsilon_2+c_6\varepsilon_{2}^{1-2(m+d)/3\delta}-c_3\leq 0$; for such a $\varepsilon_2$, and taking in mind \eqref{Auxiliar_function}, inequality \eqref{IneqPhiCasiFinal} is simplified to 
\begin{equation}\label{(eq:2.28)}
W'(t)\leq \mathcal{M}\Big(\int_\Omega V^m d{\bf x}\Big)^{\alpha}+\mathcal{N}\Big(\int_\Omega V^m d{\bf x}\Big)^{\alpha\beta}=\mathcal{M}W^\alpha+\mathcal{N}W^{\alpha\beta}\quad \textrm{on} \quad (0,t^*).
\end{equation}
Since we are assuming that $W(t)\nearrow \infty$ as $t\searrow t^{*}$, $W(t)$ can be non decreasing, so that $W(t) \geq W(0)>0$ with $t\in [0,t^{*})$, or non increasing (possibly presenting oscillations), so that there exists a time $t_1$ where $W(t_1)=W(0)$. In any case, we can write  $W(t) \geq W(0)$ for all  $t \in [t_1, t^{*})$, where $0\leq t_1<t^{*}$. By virtue of \eqref{DefinitionCoeff} of \ref{Assumption0},  $\alpha,\beta>1$, so that this implies that
\begin{equation*}
 W(t) \leq W (0)^{1-\beta} W(t)^{\beta}, \quad t \in [t_1,t^{*}),
\end{equation*} 
which, in conjunction with \eqref{(eq:2.28)}, produces
\begin{equation}\label{Diri_Ine_A}
W'(t) \leq  (\mathcal{M}W(0)^{(1-\beta)\alpha}+\mathcal{N})W^{\alpha\beta}, \quad t \,\in [t_1,t^{*}).
\end{equation}
Finally, integrating \eqref{Diri_Ine_A} between $t_1$ and $t^{*}$, we arrive at (recall $W(t_1)=W(0)$) the inequality
\begin{equation*}\label{Lower_Dirichlet}
\begin{split}
\frac{W(0)^{-\alpha\beta+1}}{-\alpha\beta+1}&=\frac{W(\tau)^{-\alpha\beta+1}}{-\alpha\beta+1}\bigg\lvert^ {t^{*}}_{t_1} \leq \int_{t_1}^{t^{*}}  (\mathcal{M}W(0)^{(1-\beta)\alpha}+\mathcal{N})d\tau\\ & \leq \int_0^{t^{*}}  (\mathcal{M}W(0)^{(1-\beta)\alpha}+\mathcal{N})d\tau
= (\mathcal{M}W(0)^{(1-\beta)\alpha}+\mathcal{N})t^{*},
\end{split}
\end{equation*}
which concludes the proof.
\end{proof}
\end{theorem}
In the behalf of scientific completeness, we point out that the previous Theorem is an extension of the main result derived in \cite{ShaeferWithoutGradient} (Schaefer, 2008), where the gradient nonlinearity for $g$ does not take part  ($k_2=0$). In this sense, in the proof of Theorem \ref{TheoremEstimateBlowUp} we use some ideas of \cite{ShaeferWithoutGradient}, but other derivations are necessarily required exactly due to the presence of $|\nabla u|^q$; these further computations imply inter alia the largeness of $k_2$ (relation \eqref{AssumptionCoeffieintTheoremBlow-Up}), not appearing in the more recent contribution.

\subsubsection*{Conflicts of interest}
The authors declare that they have no conflicts of interest.
\subsubsection*{Acknowledgements}
GV is member of the Gruppo Nazionale per l'Analisi Matematica, la Probabilit\`a e le loro Applicazioni (GNAMPA) of the Istituto Na\-zio\-na\-le di Alta Matematica (INdAM). The research of TL is supported by NNSF of P.R. China (Grant No. 61503171), CPSF (Grant No. 2015M582091), NSF of Shandong Province (Grant No. ZR2016JL021),
DSRF of Linyi University (Grant No. LYDX2015BS001), and the AMEP of Linyi University, P.R. China.


\begin{thebibliography}{10}

\bibitem{AiTsEdYaMi}
M.~Aida, T.~Tsujikawa, M.~Efendiev, A.~Yagi, and M.~Mimura.
\newblock Lower estimate of the attractor dimension for a chemotaxis growth
  system.
\newblock {\em J. London. Math. Soc.}, 74(02):453--474, 2006.

\bibitem{AndrMazSimoToledo}
F.~Andreu, J.~M. Maz\'on, F.~Simondon, and J.~Toledo.
\newblock Global existence for a degenerate nonlinear diffusion problem with
  nonlinear gradient term and source.
\newblock {\em Math. Ann.}, 314(4):703--728, 1999.

\bibitem{AndreuEtAl}
F.~Andreu, J.~M. Maz\'on, F.~Simondon, and J.~Toledo.
\newblock Blow-up for a class of nonlinear parabolic problems.
\newblock {\em Asymptot. Anal.}, 29(2):143--155, 2002.

\bibitem{ARONSON197833}
D.~Aronson and H.~Weinberger.
\newblock Multidimensional nonlinear diffusion arising in population genetics.
\newblock {\em Adv. Math.}, 30(1):33 -- 76, 1978.

\bibitem{Aronson1986}
D.~G. Aronson.
\newblock {\em The porous medium equation}, pages 1--46.
\newblock Springer Berlin Heidelberg, Berlin, Heidelberg, 1986.

\bibitem{Ball77remarkson}
J.~M. Ball.
\newblock Remarks on blow-up and nonexistence theorems for nonlinear evolution
  equations.
\newblock {\em Quart. J. Math. Oxford}, 28:473--486, 1977.

\bibitem{BandleBrunnerSurvey}
C.~Bandle and H.~Brunner.
\newblock Blowup in diffusion equations: A survey.
\newblock {\em J. Comput. Appl. Math.}, 97(1--2):3 -- 22, 1998.

\bibitem{CaoZheng}
X.~Cao and S.~Zheng.
\newblock {Boundedness of solutions to a quasilinear parabolic-elliptic
  Keller-Segel system with logistic source}.
\newblock {\em Math. Meth. Appl. Sci.}, 37(15):2326--2330, 2014.

\bibitem{Fujita_1966}
H.~Fujita.
\newblock {On the blowing up of solutions of the Cauchy problem for
  {$u_t=\Delta u+u^{1+\alpha}$}}.
\newblock {\em J. Fac. Sci. Univ. Tokyo}, 13:109--124, 1966.

\bibitem{Galaktionov1981_Russian}
V.~A. Galaktionov.
\newblock A boundary value problem for the nonlinear parabolic equation
  {$u_{t}=\Delta u^{\sigma +1}+u^{\beta }$}.
\newblock {\em Differentsial'nye Uravneniya}, 17(5):836--842, 956, 1981.

\bibitem{galaktionov_1994}
V.~A. Galaktionov.
\newblock {Blow-up for quasilinear heat equations with critical Fujita's
  exponents}.
\newblock {\em Proceedings of the Royal Society of Edinburgh: Section A
  Mathematics}, 124(3):517--525, 1994.

\bibitem{GalaktionovKurdyMikha}
V.~A. {Galaktionov}, S.~P. {Kurdyumov}, A.~P. {Mikhailov}, and A.~A.
  {Samarskii}.
\newblock {Unbounded solutions of the Cauchy problem for the parabolic equation
  {$u_t=\nabla (u^\sigma\nabla u)+u^\beta$}}.
\newblock {\em Dokl. Phys.}, 25:458--459, 1980.

\bibitem{GalaktionovShmarevVazquez}
V.~A. Galaktionov, S.~I. Shmarev, and J.~L. V\'azquez.
\newblock Behaviour of interfaces in a diffusion-absorption equation with
  critical exponents.
\newblock {\em Interfaces Free Bound.}, 2(4):425--448, 2000.

\bibitem{GalaktionovVazquezExtinction}
V.~A. Galaktionov and J.~L. V\'azquez.
\newblock Extinction for a quasilinear heat equation with absorprtion i.
  technique of intersection comparison.
\newblock {\em Comm. Partial Differential Equations}, 19(7-8):1075--1106, 1994.

\bibitem{grant2014theoryODE}
C.~Grant.
\newblock {\em Theory of Ordinary Differential Equations}.
\newblock CreateSpace Independent Publishing Platform, 2014.

\bibitem{GURTIN197735}
M.~E. Gurtin and R.~C. MacCamy.
\newblock On the diffusion of biological populations.
\newblock {\em Math. Biosc.}, 33(1--2):35 -- 49, 1977.

\bibitem{Kielhofer1974}
H.~Kielh{\"o}fer.
\newblock {Halbgruppen und semilineare Anfangs-Randwertprobleme}.
\newblock {\em Manuscripta Math.}, 12(2):121--152, 1974.

\bibitem{kobayashi1977}
K.~Kobayashi, T.~Sirao, and H.~Tanaka.
\newblock On the growing up problem for semilinear heat equations.
\newblock {\em J. Math. Soc. Japan}, 29(3):407--424, 07 1977.

\bibitem{Krylov}
N.~V. Krylov.
\newblock {\em Nonlinear elliptic and parabolic equations of the second order},
  volume~7 of {\em Mathematics and its Applications (Soviet Series)}.
\newblock D. Reidel Publishing Co., Dordrecht, 1987.

\bibitem{LSUBookInequality}
O.~A. Lady\v{z}enskaja, V.~A. Solonnikov, and N.~N. Ural'ceva.
\newblock {Linear and Quasi-Linear Equations of Parabolic Type}.
\newblock In {\em {Translations of Mathematical Monographs}}, volume~23.
  American Mathematical Society, 1988.

\bibitem{LevineTheRoleOf}
H.~A. Levine.
\newblock The role of critical exponents in blowup theorems.
\newblock {\em SIAM Rev.}, 32(2):262--288, 1990.

\bibitem{MarrasViglialoroVernier_Kodai}
M.~Marras, S.~Piro, and G.~Viglialoro.
\newblock Lower bounds for blow-up time in a parabolic problem with a gradient
  term under various boundary conditions.
\newblock {\em Kodai Math. J.}, 37(3):532--543, 2014.

\bibitem{PaynePhilVernierSECOND}
L.~Payne, G.~Philippin, and S.~V. Piro.
\newblock {Blow-up phenomena for a semilinear heat equation with nonlinear
  boundary condition, II}.
\newblock {\em Nonlinear Anal. Theory Methods Appl.}, 73(4):971--978, 2010.

\bibitem{PayneElAtl2008}
L.~Payne, G.~Philippin, and P.~Schaefer.
\newblock Blow-up phenomena for some nonlinear parabolic problems.
\newblock {\em Nonlinear Anal. Theory Methods Appl.}, 69(10):3495 -- 3502,
  2008.

\bibitem{PaynPhilSchaefer}
L.~Payne, G.~Philippin, and P.~Schaefer.
\newblock Bounds for blow-up time in nonlinear parabolic problems.
\newblock {\em J. Math. Anal. Appl.}, 338(1):438 -- 447, 2008.

\bibitem{PAYNE_Shaefer_2007}
L.~Payne and P.~Schaefer.
\newblock {Lower bounds for blow-up time in parabolic problems under Dirichlet
  conditions}.
\newblock {\em J. Math. Anal. Appl.}, 328(2):1196 -- 1205, 2007.

\bibitem{PayneSchaefer_Robin}
L.~Payne and P.~Schaefer.
\newblock {Blow-up in parabolic problems under Robin boundary conditions}.
\newblock {\em Appl. Anal.}, 87(6):699--707, 2008.

\bibitem{PaynePhiProytc}
L.~E. Payne, G.~A. Philippin, and V.~Proytcheva.
\newblock {Continuous dependence on the geometry and on the initial time for a
  class of parabolic problems I}.
\newblock {\em Math. Meth. Appl. Sci.}, 30(15):1885--1898, 2007.

\bibitem{PhilippProyt}
G.~A. Philippin and V.~Proytcheva.
\newblock Some remarks on the asymptotic behaviour of the solutions of a class
  of parabolic problems.
\newblock {\em Math. Meth. Appl. Sci.}, 29(3):297--307, 2006.

\bibitem{ProtMurrWeinberger}
M.~H. Protter and H.~F. Weinberger.
\newblock {\em Maximum principles in differential equations}.
\newblock Springer-Verlag, New York, 1984.

\bibitem{Sacks}
P.~E. Sacks.
\newblock Global behavior for a class of nonlinear evolution equations.
\newblock {\em SIAM J. Math. Anal.}, 16(2):233--250, 1985.

\bibitem{ShaeferWithoutGradient}
P.~Schaefer.
\newblock Lower bounds for blow-up time in some porous medium problems.
\newblock {\em Proc. Dynam. Systems Appl.}, 5:442--445, 2008.

\bibitem{ShaeferExistClassicaPorous}
P.~Schaefer.
\newblock Blow-up phenomena in some porous medium problems.
\newblock {\em Dynam. Systems Appl.}, 18:103--110, 2009.

\bibitem{Souplet_Gradient}
P.~Souplet.
\newblock Finite time blow-up for a non-linear parabolic equation with a
  gradient term and applications.
\newblock {\em Math. Meth. Appl. Sci.}, 19(16):1317--1333, 1996.

\bibitem{Talenti1976BestConstant}
G.~Talenti.
\newblock {Best constant in Sobolev inequality}.
\newblock {\em Annali di Matematica Pura ed Applicata}, 110(1):353--372, 1976.

\bibitem{vazquez2007porous}
J.~V\'azquez.
\newblock {\em The Porous Medium Equation: Mathematical Theory}.
\newblock Oxford Mathematical Monographs. Clarendon Press, 2007.

\bibitem{Verhulst}
P.-F. Verhulst.
\newblock Notice sur la loi que la population poursuit dans son accroissement.
\newblock {\em Correspondance math{\'e}matique et physique}, 10:113--121, 1838.

\bibitem{ViglialoroDifferentialIntegralEquations}
G.~Viglialoro.
\newblock {Blow-up time of A Keller-Segel-type system with Neumann and Robin
  boundary conditions}.
\newblock {\em Differ. Integral Equ.}, 29(3-4):359--376, 2016.

\bibitem{ViglialoroBoundnessVeryWeak}
G.~Viglialoro.
\newblock Boundedness properties of very weak solutions to a fully parabolic
  chemotaxis-system with logistic source.
\newblock {\em Nonlinear Anal. Real World Appl.}, 34:520--535, 2017.

\bibitem{ViglialoroWolleyDCDS}
G.~Viglialoro and T.~Woolley.
\newblock Eventual smoothness and asymptotic behaviour of solutions to a
  chemotaxis system perturbed by a logistic growth.
\newblock {\em Discrete Continuous Dyn. Syst. Ser. B.}, 22(5), 2017.

\bibitem{WinklerVeryWeakSol}
M.~Winkler.
\newblock Chemotaxis with logistic source: very weak global solutions and their
  boundedness properties.
\newblock {\em J. Math. Anal. Appl.}, 348(2):708--729, 2008.

\bibitem{WinDespiteLogistic}
M.~Winkler.
\newblock {Blow-up in a higher-dimensional chemotaxis system despite logistic
  growth restriction}.
\newblock {\em J. Math. Anal. Appl.}, 384(2):261--272, 2011.

\end{thebibliography}
\end{document}